\documentclass[leqno]{amsart}
\usepackage[cp1250]{inputenc}
\usepackage[T1]{fontenc}
\usepackage{amsthm}
\usepackage[leqno]{amsmath}
\usepackage{amssymb}
\usepackage{amsfonts}
\usepackage{enumerate}
\usepackage{color}

\newcommand{\kk}{\mathbf{\mathbb{K}}}
\newcommand{\zz}{\mathbf{\mathbb{Z}}}
\newcommand{\rr}{\mathbf{\mathbb{R}}}
\newcommand{\cc}{\mathbf{\mathbb{C}}}

\newcommand{\vf}{\varphi}

\newcommand{\dist}{\operatorname{dist}}

\newtheorem{lemma}{Lemma}
\newtheorem{theorem}{Theorem}

\newtheorem{cor}{Corollary}

\theoremstyle{remark}
\newtheorem{definition}{Definition}
 
\title{Sufficiency of non-isolated singularities}

\author{Piotr Migus, Tomasz Rodak, Stanis{\l}aw Spodzieja}
\thanks{This research was partially supported by the Polish OPUS Grant No 2012/07/B/ST1/03293}

\date{\today}

\begin{document}

\begin{abstract}
  We give, in terms of the {\L}ojasiewicz inequality, a sufficient condition for $C^k$-mappings germs of non-isolated singularity at zero to be isotopical. 
\end{abstract}
\maketitle

\section{Introduction and results}

Let $F:(\rr^n,a)\to \rr^m$ denote a mapping defined in a neighbourhood of $a\in \rr^n$ with values in $\rr^m$. If  $F(a)=b$, we put $F:(\rr^n,a)\to (\rr^m,b)$. By $\nabla f$ we denote the gradient of a $C^1$-function $f:(\rr^n,a)\to \rr$. By $|\,\cdot\,|$ we denote a norm in $\rr^n$ and $\dist(x,V)$ - the distance of a point $x\in \rr^n$ to a set $V\subset \rr^n$ (or $\dist (x,V)=1$ if $V=\emptyset$).


By a $k$-\emph{jet} at $a\in \rr^n$ in the $C^l$ class we mean a family of $C^l$ functions $(\rr^n,a)\to \rr$, called \emph{$C^l$-realisations} of this jet, possesing the same Taylor polynomial of degree $k$ at $a$. The $k$-jet is said to be $C^r$-\emph{sufficient} (respectively $C^r$-$v$-\emph{sufficient}) in the $C^l$ class, if for every of his $C^l$-realisations $f$ and $g$ there exists a $C^r$ diffeomorphism $\vf:(\rr^n,a)\to(\rr^n,a)$, such that $f\circ \vf=g$ (respectively $f^{-1}(0)=\vf(g^{-1}(0))$\,) in a neighbourhood of $a$ (R. Thom \cite{Th}). 

In the paper we will consider the $k$-jets in the class $C^k$ and write shortly - $k$-jets.

The clasical result in the subject sufficiency of jets is the following:

\begin{theorem}[Kuiper, Kuo, Bochnak-{\L}ojasiewicz]\label{KKBL}
Let $w$ be a $k$-jet at $0\in \rr^n$ and let $f$ be its $C^k$ realisation. If $f(0)=0$  then the following conditions are equivalent:
\begin{itemize}
\item[(a)] $w$ is $C^0$-sufficient in the $C^k$ class,
\item[(b)] $w$ is $C^0$-$v$-suffucient in the $C^k$ class,
\item[(c)] $|\nabla f(x)|\ge C|x|^{k-1}$ as $x\to 0$ for some constant $C>0$.
\end{itemize}
\end{theorem}

The implication (c)$\Rightarrow$(a) was proved by N.~H.~Kuiper \cite{Kui} and T.~C.~Kuo \cite{Kuo}, (b)$\Rightarrow$(c) - by J.~Bochnak and S. \L ojasiewicz \cite{BL}, and the implication (a)$\Rightarrow$(b) is obvious (see also \cite{OSS}, \cite{Sp1}). Analogous result in the complex case was proved by S. H. Chang and Y. C. Lu \cite{CL}, B. Teissier \cite{Te} and J. Bochnak and W.~Kucharz \cite{BK}. Similar considerations as above are carried out for functions in a neighbourhood of infinity (see \cite{CH}, \cite{Sk}, \cite{RS}). 


Theorem \ref{KKBL} concerns the isolated singularity of $f$ at $0$, i.e. the point $0$ is an isolated zero of $\nabla f$. The case of non-isolated singularities of real functions was investigated by many authors, for instance by J. Damon and T. Gaffney \cite{DG},  T.~~Fukui and E.~~Yoshinaga \cite{FY}, V. Grandjean \cite{G}, Xu Xu \cite{Xu} and for complex functions - by D.~~Siersma \cite{Si1,Si2} and R. Pellikaan \cite{P}. 

The purposes of this article are generalisations of the above results for a $C^k$ mappings in a neighbourhood of zero with non-isolated singularity at zero. Recall the definition of $k$-$Z$-jet in the class of functions with non-isolated singularity at zero (cf. \cite{Xu}). 

The set of $C^k$ mappings $(\rr^n,a)\to \rr^m$ we denote by $\mathcal{C}^k_a(n,m)$. For a function $f\in \mathcal{C}^k_a(n,1)$, by $j^kf(a)$ we denote the $k$-jet at $a$ (in the $C^k$-class) determined by $f$. For a mapping $F=(f_1,\ldots,f_m)\in \mathcal{C}^k_a(n,m)$ we put $j^kF(a)=(j^kf_1(a),\ldots,j^kf_m(a))$.

Let $Z\subset \rr^n$ be a set such that $0\in Z$ and let $k\in\zz$, $k>0$. By $k$-$Z$-\emph{jet in the class $\mathcal{C}^k_0(n,m)$}, or shortly $k$-$Z$-\emph{jet}, we mean an equivalence class $w\subset \mathcal{C}^k_0(n,m)$ of the equivalence relation $\sim$: $F\sim G$ iff for some neighbourhood $U\subset \rr^n$ of the origin, $j^kF(a)=j^kG(a)$ for $a\in Z\cap U$ (cf. \cite{Xu}). The mappings $F\in w$ we call \emph{$C^k$-$Z$-realisations of the jet} $w$ and we write $w=j^{k}_{Z}F$. The set of all jets $j^{k}_{Z}F$ we denote by $J^{k}_Z(n,m)$. 

The $k$-$Z$-jet $w\in J^{k}_Z(n,m)$ is said to be $C^r$-$Z$-\emph{sufficient} (resp. $C^r$-$Z$-$v$-\emph{sufficient}) in the $C^k$ class, if for every of its $C^k$-$Z$-realisations $f$ and $g$ there exists a $C^r$ diffeomorphism $\vf:(\rr^n,0)\to(\rr^n,0)$, such that $f\circ \vf=g$ (resp. $f^{-1}(0)=\vf(g^{-1}(0))$\,) in a neighbourhood $U$ of $0$ and $\vf(x)=x$ for $x\in Z\cap U$.

The following Kuiper and Kuo criterion (Theorem \ref{KKBL} (c)$\Rightarrow$(a)) for jets with non-isolated singularity was proved by Xu Xu \cite{Xu}.

\begin{theorem}\label{XuXu}
Let $Z\subset \mathbb{R}^n$ be a 
closed set 
 such that $0\in Z$. 
  If $f\in \mathcal{C}^k(n,1)$ such that $\nabla f(x)=0$ for $x\in Z$, satisfies the condition 
\begin{equation}\label{KKcondition}
|\nabla f(x)|\ge C\dist(x,Z)^{k-1}\hbox{ as $x\to 0$ for some constant $C>0$},
\end{equation}
then the $k$-$Z$-jet of f is $C^0$-$Z$-sufficient.
\end{theorem}

The main result of this paper is Theorem \ref{main.result} below. It is a generalisation of the Theorem \ref{XuXu} to the case of mapping jets. Let us start with some definition. Let $X,Y$ be Banach spaces over $\rr$. Let ${L}(X,Y)$ denote the Banach space
of linear continuous mappings from $X$ to $Y$. For $A\in {L}(X,Y)$, $A^*$ stands for the
adjoint operator in ${L}(Y',X')$, where $X'$ is the dual space of $X$. For $A\in{L}(X,Y)$ we put
  \begin{equation}
    \label{rabier.function}
    \nu(A)=\inf\{\|A^*\vf\|:\vf\in Y',\|\vf\|=1\},
  \end{equation}
  where $\|A\|$ is the norm of linear mapping $A$ (see \cite{Ra}). In the case $f\in\mathcal{C}^k_0(n,1)$ we have $\nu (d f)=|\nabla f|$, where $d f $ is the differental of $f$.


\begin{theorem}\label{main.result}
  Let $f\colon(\rr^n,0)\to(\rr^m,0)$, where $m\leq n$, be a $C^k$-$Z$-realisation of a $k$-$Z$-jet $w\in J^{k}_Z(n,m)$, where $k>1$ and $Z=\{x\in \rr^n:\nu (df (x))=0\}$, $0\in Z$. Assume that for a positive constant $C$, 
  \begin{equation}\label{main.assumption}
    \nu(df(x))\geq C\dist(x,Z)^{k-1} \quad \hbox{as $x\to0$}.
  \end{equation}
  Then the jet $w$ is $C^0$-$Z$-sufficient in the class $C^k$. Moreover for any $C^k$-$Z$-realisations $f_1,f_2$ of $w$,  the deformation $f_1+t(f_2-f_1)$, $t\in\rr$ is topologically trivial along $[0,1]$. In particular the mappings $f_1$ and $f_2$ are  isotopical at zero. 
\end{theorem}

For the definition of isotopy and topological triviality see Subsection \ref{Isotopysection}. By Lemmas \ref{dist.lemma} and \ref{real.structure} in Section \ref{proofofmain.result}, Theorem \ref{main.result} is also true for holomorphic mappings. It is not clear to the authors if the inverse to Theorem \ref{main.result} holds.  In the proof of Theorem~\ref{main.result}, given in Section \ref{proofofmain.result}, we use a method of the proof of Theorem~~1 in \cite{RS}.

In the case of nondegenerate analytic functions $f$, $g$, a conditions for topological triviality of deformations $f+tg$, $t\in[0,1]$ in terms of Newton polyhedra was obtained by J. Damon and T. Gaffney \cite{DG}, and for blow analytic triviality -- T. Fukui and E.~~Yoshinaga \cite{FY} (see also \cite{Te0}, \cite{Y}).

From the proof of Theorem \ref{main.result} we obtain a version of the theorem for functions of $C^1$ class with locally Lipschitz differentials.

\begin{cor}\label{corc1lipschita}
Let  $f,f_1:(\rr^n,0)\to(\rr^m,0)$ be differentiable mappings with locally Lipschitz differentials  $df,\,df_1:(\rr^n,0)\to {L}(\rr^n,\rr^m)$ let $Z=\{x\in \rr^n:\nu (df (x))=0\}$, and let $0\in Z$. It 
 \begin{equation}\label{assumption0}
    \nu(df(x))\geq C\dist(x,Z),
 \end{equation}
\begin{equation}\label{assumption1}
|f(x)-f_1(x)|\leq C_1\nu(df(x))^2,
\end{equation}
\begin{equation}\label{assumption2}
     \|df(x)-df_1(x)\|\leq C_2\nu(df(x))
\end{equation}
as $x\to 0$ for some constants 
 $C,C_1,C_2>0$, $C_2<\frac{1}{2}$, 
  then the deformation $f+t(f_1-f)$ is topologically trivial along $[0,1]$. In particular $f$ and $f_1$ are isotopical at zero.
\end{cor}

The proof of the above corollary is given in Subsection \ref{proofofcorollarylipschitz}.

In Section \ref{vsufficiency} we prove the following theorem type of Bochnak-{\L}o\-ja\-sie\-wicz (cf. implication (b)$\Rightarrow$(c) in Theorem \ref{KKBL}), that $C^0$-$Z$-$v$-sufficiency of a jets implies the {\L}ojasiewicz inequality, provided $j^{k-1}f(0)=0$ for $C^k$-$Z$-realisations $f$ of the jet. 
Namely, we will prove the following
   
\begin{theorem}\label{BLplus1}
Let $Z\subset \mathbb{R}^n$ be a set such that $0\in Z$, let $w$ be a $k$-$Z$-jet, $k>1$, and let $f$ be its ${C}^k$-$Z$-realisation. If  $w$ is ${C}^0$-$Z$-$v$-sufficient in ${C}^k$-class, $j^{k-1}f(0)=0$ 
 and $V(\nabla f)\subset Z$, then
\begin{equation}\label{eqBLplus}
|\nabla f(x)|\ge C\dist(x,Z)^{k-1}\hbox{ as $x\to 0$ for some constant $C>0$}.
\end{equation}
\end{theorem}

It is obvious that a $C^0$-$Z$-sufficient jet is also a $C^0$-$Z$-$v$-sufficient, so, Theorem~~\ref{BLplus1}, in a certain sense is an inverse of Theorem \ref{XuXu}. 

\section{Proof of Theorem \ref{main.result}}\label{proofofmain.result}

\subsection{Differential equations} Let us start from recalling the following

\begin{lemma}\label{lemmauniqueness}
Let $G \subset \mathbb{R}\times \mathbb{R}^n$ be an open set, $W:G\rightarrow \mathbb{R}^n$ be a continuous mapping and let $V \subset \mathbb{R}^n$ be a closed set. If in $G\backslash (\mathbb{R}\times V)$ system
\begin{equation}\label{eqlemmaunique}
\frac{dy}{dt}=W(t,y)
\end{equation}
has a global unique solutions and there exist neighbourhood $U \subset G$ of set $(\mathbb{R}\times V)\cap G$ and a positive constant $C$ such that 
\begin{equation}
\left|W(t,x) \right|\leq C \dist(x,V)\quad \textrm{ for } (t,x)\in U,
\end{equation}
then the system \eqref{eqlemmaunique} in $G$ has a global unique solutions.
\end{lemma}

\subsection{The Rabier function} 
Let $X,Y$ be Banach spaces over $\kk$, where $\kk=\rr$ or $\kk =\cc$. Let ${L}(X,Y)$ denote the Banach space
of linear continuous mappings from $X$ to $Y$. For $A\in {L}(X,Y)$, $A^*$ stands for the
adjoint operator in ${L}(Y',X')$, where $X'$ is the dual space of $X$. We begin with recalling some properties of the Rabier function (cf. \cite{RS}).

\begin{lemma}[\cite{KOS}]
  \label{dist.lemma}
  Let $\Sigma$ be the set of operators $A\in {L}(X,Y)$ such that $A(X)\subsetneq Y$.
  We have 
  \begin{equation*}
    \nu(A)=\operatorname{dist}(A,\Sigma),\quad A\in {L}(X,Y).
  \end{equation*}
\end{lemma}

\begin{lemma}[\cite{Ra}]\label{nuLipschitz} 
Let $A,B\in {L}(X,Y)$. Then $|\nu(A)-\nu(B)|\le \|A-B\|$. In particular $\nu:{L}(X,Y)\to \rr$ is Lipschitz.
\end{lemma}

From Lemma \ref{nuLipschitz} we have
\begin{lemma}
  \label{obvious.inequality}
  If $A,B\in{L}(X,Y)$ then
  \begin{equation*}
    \nu(A+B)\geq\nu(A)-\|B\|.
  \end{equation*}
\end{lemma}

\begin{definition}[\cite{J}]\label{gaffney}
  Let $\textbf{a}=[a_{ij}]$ be the matrix of $A\in{L}(\kk^n,\kk^m)$, $n\geq m$. By $M_I(A)$, 
  where $I=(i_1,\ldots,i_m)$ is any subsequence of $(1,\ldots,n)$, we denote an $m\times m$ minor of
  $\textbf a$ given by columns indexed by $I$. Moreover, if $J=(j_1,\ldots,j_{m-1})$ is any
  subsequence of $(1,\ldots,n)$ and $j\in\{1,\ldots,m\}$, then by $M_J(j)(A)$ we denote 
  an $(m-1)\times(m-1)$ minor of $\textbf{a}$ given by columns indexed by $J$ and with deleted $j$th row
  (if $m=1$ we put $M_J(j)(A)=1$). Let
  \begin{gather*}
    h_I(A)=\max\left\{|M_J(j)(A)|:J\subset I,j=1,\ldots,m\right\},\\
    g'(A)=\max_I\frac{|M_I(A)|}{h_I(A)}.
  \end{gather*}
  Here we put $0/0=0$. If $m=n$, we put $h_I=h$.
\end{definition}

\begin{lemma}[\cite{J}]\label{jelonek.lemma}
  There exist $C_1,C_2>0$, such that for any $A\in{L}(\kk^n,\kk^m)$ we have
  \begin{equation*}
    C_1g'(A)\leq \nu(A)\leq C_2g' (A).
  \end{equation*}
\end{lemma}

\begin{cor}[\cite{RS}]\label{continuity} 
  The function $g'$ is continuous.
\end{cor}

\begin{lemma}[\cite{KOS}]\label{real.structure}
  Assume that $X,Y$ are complex Banach spaces. 
  Let $\Sigma_{\cc}$ (resp. $\Sigma_\rr$) be the set of nonsurjective
  $\cc$-linear (resp. $\rr$-linear) continuous maps from $X$ to $Y$.  Then for any continuous
  $\cc$-linear map $A\colon X\to Y$,
  \begin{equation*}
    \operatorname{dist}(A,\Sigma_\cc)=\operatorname{dist}(A,\Sigma_\rr).
  \end{equation*}
\end{lemma}

\subsection{Isotopy and triviality}\label{Isotopysection}  
  Let  $\Omega\subset\rr^n$ be a neighbourhood of $0\in\rr^n$ and let $Z\subset \rr^n$ be a set such that $0\in Z$. 
  
We will say, that a
  continuous mapping $H\colon\Omega\times [0,1]\to\rr^n$ is an \emph{isotopy near $Z$ at zero} if
  \begin{enumerate}[(a)]
  \item $H_0(x)=x$ for $x\in \Omega$ and $H_t(x)=x$ for $t\in[0,1]$ and $x\in \Omega\cap Z$,
  \item for any $t$ the mapping $H_t$ is a homeomorphism onto $H_t(\Omega)$,
  \end{enumerate}
where the mapping $H_t:\Omega\to \rr^n$ is defined by $H_t(x)=H(x,t)$ for $x\in\Omega$, $t\in [0,1]$.



Let $f\colon\Omega_1\to\rr^m$, $g\colon\Omega_2\to \rr^m$ where 
  $\Omega_1,\Omega_2\subset\rr^n$ are neighbourhoods of $0\in\rr^n$ and let $Z\subset \rr^n$ be a set such that $0\in Z$. We call $f$ and $g$ \emph{isotopical near $Z$ at zero} if there exists an isotopy near $Z$ at zero 
  $H:\Omega\times [0,1]\to\rr^n$, $\Omega\subset\Omega_1\cap \Omega_2$, such that
  $f(H_1(x))=g(x)$, $x\in\Omega$. 

Let $h:\Omega_3\to \rr^m$, where $\Omega_3\subset \rr^n$ is a neighbourhood of $0\in\rr^n$. We say that a deformation $f+th$, is \emph{topologically trivial near $Z$} along $[0,1]$ if there exists an isotopy near $Z$ at zero $H:\Omega \times [0,1]\to \rr^m$, $\Omega \subset \Omega_1\cap\Omega_2$, such that $f(H(t,x))+th(H(t,x))$ do not depend on $t$.

%
%

\subsection{Proof of Theorem \ref{main.result}} By $dP$ we denote the differential of $P$ and $d P(x)$ -- the diferential of $P$ at the point $x$. By $d_x P$ we denote the differential of $P$ with respect to the system of variables $x$. 
 

Let $f,f_1\in w$ and let $P=f_1-f=(P_1,\ldots,P_m)$. Then we have $j^kP(a)=0$ for $a\in Z\cap U$ for some neighbourhood $U\subset \rr^n$ of $0$. In consequence, decreasing if necessary $U$, we may assume that 
  \begin{equation}\label{eqP1}
  |P(x)|\le \frac {C}{3}\dist (x,Z)^k\quad\hbox{and}\quad\|dP(x)\|\le \frac{C}{3}\dist(x,Z)^{k-1}
  \end{equation}
for $x\in U$.

  Consider the mapping $F\colon\rr\times\rr^n\to\rr^m$,
  \begin{equation*}
    F(\xi,x)=f(x)+\xi P(x).
  \end{equation*}
  Let us fix $\xi\in (-2,2)$. By \eqref{eqP1} and Lemma \ref{obvious.inequality} we get
  \begin{equation*}\label{eqforcorollary1}
      \nu(d_xF(\xi,x))
      \geq\nu(df(x))-|\xi|\|dP(x)\|\geq\frac{C}{3}\dist(x,Z)^{k-1},\quad x \in U.
  \end{equation*}
  Thus by Lemma \ref{jelonek.lemma} there exists $C'>0$ such that
  \begin{equation}
    \label{dxF.inequality}
        g'(d_xF(\xi,x))\geq C' \dist(x,Z)^{k-1},\quad \xi\in(-2,2),\quad x\in U.
  \end{equation}
  Set $G=\{(\xi,x)\in\rr\times U:|\xi|<2 \}$. 
  In the notation of Definition \ref{gaffney} we put 
  \begin{equation*}
    \begin{split}
      A_I=\left\{(\xi,x)\in G:\frac{|M_I(d_xF(\xi,x))|}{h_I(d_xF(\xi,x))}\leq\frac{C'}{2}\dist(x,Z)^{k-1}\right\}.
    \end{split}
  \end{equation*}
  By Corollary \ref{continuity}
  the sets $A_I$ are closed in $G$ and $(\rr\times Z)\cap G\subset A_I$. From  \eqref{dxF.inequality} 
  we see that $\{G\setminus A_I:I\}$ is an open covering of $G\setminus (\rr\times Z)$. Let
  $\{\delta_I:I\}$ be a $C^\infty$ partition of unity associated to this covering.
  
  Let us consider the following system of linear equations 
  \begin{equation}\label{linear.system}
    (d_xF(\xi,x))W(\xi,x)^T=-P(x)^T.
  \end{equation}
with indeterminates $W(\xi,x)=(W_1(\xi,x),\ldots,W_n(\xi,x))$ and parameters $(\xi,x)\in G$. 
  Let us take any subsequence $I=(i_1,\ldots,i_m)$ of the sequence $(1,\ldots,n)$. 
  For simplicity of notation we assume that $I=(1,\ldots,m)$. For all $(\xi,x)\in G$
  such that $M_I(d_xF(\xi,x))\ne 0$ we put 
  $W^I(\xi,x)=(W^I_1(\xi,x),\ldots,W^I_n(\xi,x))$, by
  \begin{equation*}
    \begin{split}
      W^{I}_l(\xi,x)&=\sum_{j=1}^m(-P_j(x))(-1)^{l+j}\frac{M_{I\setminus l}(j)(d_xF(\xi,x))}{M_I(d_xF(\xi,x))},
      \quad l=1,\ldots,m,\\
      W^I_l(\xi,x)&=0,\quad l=m+1,\ldots,n,
    \end{split}
  \end{equation*}
where $I\setminus l=(1,\ldots,l-1,l+1,\ldots,m)$ for $l=1,\ldots,m$.  
  Cramer's rule implies 
  \begin{equation*}
    (d_xF(\xi,x))W^I(\xi,x)^T=-P(x)^T.     
  \end{equation*}
Since $k>1$, then 
 $\delta_IW^I$ is a $C^1$ mapping on $G\setminus (\rr\times Z)$ (after suitable extension). Hence 
  $W=\sum_I\delta_IW^I$ is also $C^1$ mapping 
  on $G\setminus(\rr\times  Z)$. We put $W(\xi,x)=0$ for $(\xi,x)\in (\rr \times Z)\cap G$.  It is easy to see, that 
  $W$ satisfies the equation \eqref{linear.system}. 

Observe that 
  \begin{equation}
    \label{W.inequality}
    \|W(\xi,x)\|\leq C''\dist(x,Z),\quad \xi\in(-2,2),\quad x\in U,
  \end{equation}
  where $C''=2mC\sqrt{n}/(3C')$. Indeed from \eqref{eqP1} the definitions of $A_I$, the choice of $P$ and the above construction we get 
  \begin{equation*}
    \begin{split}
      \|W(\xi,x)\|
      \leq&\sum_{\{I:\delta_I(\xi,x)\ne 0\}}\delta_I(\xi,x)\|W^I(\xi,x)\|\\
      \leq& \sum_{\{I:\delta_I(\xi,x)\ne 0\}}\delta_I(\xi,x) \sqrt{n}\max_{l=1}^m|W_l^I(\xi,x)|\\
      \leq& \sum_{\{I:\delta_I(\xi,x)\ne 0\}}\delta_I(\xi,x)  \sqrt{n} \sum_{j\in I}|P_j(x)|\frac{h_I(d_xF(\xi,x))}{|M_I(d_xF(\xi,x))|}\\
      \leq& \sum_{\{I:\delta_I(\xi,x)\ne 0\}}\delta_I(\xi,x)  \sqrt{n} \sum_{j\in I}\frac{C}{3} \dist(x,Z)^k\frac{2}{C'}\frac{1}{\dist(x,Z)^{k-1}}\\
      =\, & m\sqrt{n}\frac{C}{3}\frac{2}{C'}\dist(x,Z).
    \end{split}
  \end{equation*}

  Let us consider the following system of differential equations
  \begin{equation}
    \label{differential.equation}
    y'=W(t,y).
  \end{equation}
  Since $W$ is at least of class $C^1$ on $G\setminus(\rr\times Z)$, so it is a locally lipschitzian vector field. As a consequence, the above system has a uniqueness of solutions property in $G\setminus(\rr\times Z)$, Hence, inequality  \eqref{W.inequality} and Lemma \ref{lemmauniqueness} implies the global uniqueness of solutions of the system \eqref{differential.equation} in $G$. 
  
  Choose $(\xi,x)\in G$ and define $\vf_{(\xi,x)}$ to be the maximal solution of
  \eqref{differential.equation} such that $\vf_{(\xi,x)}(\xi)=x$. Set 
  $\Omega_0=\{x\in\rr^n:\|x\|<r_0\}$, $\Omega_1=\{x\in\rr^n:\|x\|<r_1\}$, where $r_0,r_1>0$. Since $0\in Z$, the mapping $\vf(\xi)=0$, $\xi\in\rr$ is a solution of \eqref{differential.equation}. Hence for sufficiently small $r_0,r_1$, for any $x\in \Omega_0$, the solution  $\vf_{(0,x)}$ is defined on $[0,1]$ and $\vf_{(0,x)}(t)\in\Omega_1$, if $t\in[0,1]$ and for any $x\in \Omega_1$, the solution $\vf_{(1,x)}$ is also defined on $[0,1]$. 
  Let $H,\widetilde H\colon \Omega_0\times[0,1]\to \Omega_1$ be given by
  \begin{equation*}
    H(x,t)=\vf_{(0,x)}(t),\quad
    \widetilde H(y,t)=\vf_{(t,y)}(0).
  \end{equation*}
  The mappings $H,\widetilde H$ are well defined. 
  Moreover one can extend these mappings to continuous mappings on some open neighbourhood of 
  $\Omega_0\times [0,1]$. Put $\Omega=\Omega_1$, 
  $\Omega^t=\{y\in\rr^n:\widetilde H(y,t)\in \Omega\}$, $t\in [0,1]$.
  By uniqueness solutions of \eqref{differential.equation} for any $t$ we have
  $\widetilde H(H(x,t),t)=x$, 
  $H(x,0)=x$,  $x\in \Omega$, and $H(\widetilde H(y,t))=y$, $y\in\Omega^t$. 
  Moreover there exists a neighbourhood $\Omega'\subset \rr^n$ of $0$ such that 
  $\Omega'\subset \Omega^t$ for any $t$.
  
  Finally, by \eqref{linear.system} we have
  \begin{equation*}
      \frac{d}{dt}F(t,\vf_{(\xi,x)}(t))^T
      =P(x)^T+(d_xF)(t,\vf_{(\xi,x)}(t))W(t,\vf_{(\xi,x)}(t))^T=0,
  \end{equation*}
  so, $F(t,\vf_{(0,x)}(t))=f(x)$ 
  and consequently 
  $f(H(x,1))+tP(H(x,1))=f(x)$ for $t\in[0,1]$ and  $x\in \Omega'$. This ends the proof.\hfill$\square$

\subsection{Proof of Corollary \ref{corc1lipschita}}\label{proofofcorollarylipschitz}
Under notations of the proof of Theorem \ref{main.result}, by \eqref{assumption0}, \eqref{assumption2} and Lemma \ref{nuLipschitz} we obtain 
     $ \nu(d_xF(\xi,x))=\nu(df(x)+\xi dP(x)) \geq\nu(df(x))-|\xi|\|dP(x)\|\geq C(1-2C_2)\dist(x,Z)$, $ x \in U$.
     Obviously $C(1-2C_2)>0$. 
Then there exists $C'>0$ such that
  \begin{equation}\label{dxF.inequality1}
        g'(d_xF(\xi,x))\geq C' \dist(x,Z),\quad \xi\in(-2,2),\quad x\in U.
  \end{equation}
So, we will use \eqref{dxF.inequality1} instead of \eqref{dxF.inequality}. By \eqref{assumption1} we obtain \eqref{W.inequality}. Moreover, the assumption that $df$ and $df_1$ are locally Lipschitz mappings implies that the mapping $W$ is locally Lipschitz outside $(-2,2)\times Z$. Then, by the same argument as in the proof of Theorem  \ref{main.result} we deduce the assertion.

\section{Proof of Theorem \ref{BLplus1}}\label{vsufficiency}

We will use the idea from \cite{BL}. It suffices to prove of the Theorem for $Z=V(\nabla f)$. Suppose to the contrary that for any neighbourhood $U$ of $0$ and for any constant $C>0$ there exist $x \in U$ such that
$$
\left|\nabla f(x) \right|< C \dist(x,Z)^{k-1}.
$$
Then for some sequence $(a_{\nu})\subset \mathbb{R}^n\backslash Z$ such that $a_{\nu}\rightarrow 0$ when $\nu \rightarrow \infty$ we have
\begin{equation}\label{eqeq1}
\left|\nabla f(a_{\nu}) \right|\le \frac{1}{\nu}\dist(a_{\nu},Z)^{k-1} \quad\textrm{ for  }\nu \in \mathbb{N}.
\end{equation}
Choosing a subsequence of $(a_{\nu})$, if necessary,  we can assume that 
$$
\dist(a_{\nu+1},Z) < \frac{1}{2}\dist(a_{\nu},Z), \quad \textrm{ for } \nu \in \mathbb{N}.
$$
Then
$$
B_{\nu}=\{x\in \mathbb{R}^n:\left|x-a_{\nu}\right|\leq \frac{1}{4}\dist (a_{\nu},Z)\}, \quad \nu \in\mathbb{N},
$$
is family of pairwise disjoint balls.

Let us take sequence $(\lambda_{\nu})\subset \mathbb{R}$ such that $\lambda_{\nu}>0$ for any $\nu\in\mathbb{N}$ and
\begin{equation}\label{eqeq3}
\frac{\lambda_{\nu}}{\dist(a_{\nu},Z)^{k-2}}\rightarrow 0,\quad \nu \rightarrow \infty.
\end{equation} 
Since $k>1$, we may assume that 
\begin{equation}\label{eqeq4}
\lambda_{\nu}\textrm{  is not eigenvalue of matrix } \left[ \frac{\partial^2 f }{\partial x_i \partial x_j}(a_{\nu})\right].
\end{equation}
Let $\alpha:\mathbb{R}^n\rightarrow \mathbb{R}$ be function of $\mathcal{C}^{\infty}$-class such that $\alpha(x)=0$ for $|x|\geq \frac{1}{4}$ and $\alpha(x)=1$ in some neigbourhood of $0$. By $\langle \cdot,\cdot\rangle$ we denote the standard inner product in $\rr^n$. Consider function $F:\mathbb{R}^n\rightarrow \mathbb{R}$ defined by the formulas
$$
F(x)=\alpha\left( \frac{x-a_{\nu}}{\dist(a_{\nu},Z)}\right)\left( f(a_{\nu})+ \langle \nabla f(a_\nu), x-a_\nu\rangle 
+\frac{1}{2}\lambda_{\nu} |x-a_{\nu}|^2\right),
$$
for $x\in B_{\nu}$ and $F(x)=0$ for $x\notin \bigcup_{\nu=1}^{\infty}B_{\nu}$. Then $F$ is a $\mathcal{C}^{k}$-function and $F(0)=0$. Moreover $f(a_{\nu})=F(a_{\nu})$ and $\nabla f(a_{\nu})=\nabla F(a_{\nu})$ so 
\begin{equation}\label{eqeq5}
(f-F)(a_{\nu})=0 \quad \textrm{ and }\quad \nabla(f-F)(a_{\nu})=0, \quad \nu \in \mathbb{N}.
\end{equation}

Let $M>0$ be such that $|\alpha(x)|\leq M$ for $x\in \mathbb{R}^n$. Then for $x\in B_{\nu}$ we have
\begin{align*}
\frac{|F(x)|}{\dist(x,Z)^k} &\leq M \frac{|f(a_{\nu})+ \langle \nabla f(a_\nu), x-a_\nu\rangle 
+\frac{1}{2}\lambda_{\nu} |x-a_{\nu}|^2|}{\dist(x,Z)^k}\\
& \leq 2^k M\frac{|f(a_{\nu})|+|\nabla f(a_{\nu})|\dist(a_{\nu},Z)+\frac{1}{2}|\lambda_{\nu}|\dist(a_{\nu},Z)^2}{\dist(a_{\nu},Z)^k}.
\end{align*}
Since $j^{k-1}f(0)=0$, 
 then 
\begin{equation*}
\frac{\left| f(a_{\nu}) \right|}{|a_{\nu}|^{k-1}}\rightarrow 0, \quad \textrm{ when  }\nu \rightarrow \infty.
\end{equation*}
Hence, from the above, and from \eqref{eqeq1} 
 and \eqref{eqeq3} we obtain
$$\frac{|F(x)|}{\dist(x,Z)^k} \rightarrow 0,\quad \textrm{ when } x\rightarrow 0,$$
so
$$
\frac{|F(x)|}{|x|^k} \rightarrow 0,\quad \textrm{ when } x\rightarrow 0.
$$
Therefore $f-F$ is ${C}^k$-$Z$-realisation of $k$-$Z$-jet $w$ (recall that for any $x\in Z\setminus\{0\}$ the function  $F$ vanishes in a neighbourhood of $x$). From \eqref{eqeq5} we have that $(f-F)$ has zeros outside the set $Z$, so by our assumption, $f$ has zeros outside the set $Z$. By the implicit function theorem for some neighbourhood $U$ of $0\in \mathbb{R}^n$, we obtain that $f^{-1}(0)\cap(U\setminus Z)$ is $(n-1)$-dimensional topological manifold. Therefore $(f-F)^{-1}(0)\cap(U_1\backslash Z)$ is also $(n-1)$-dimensional topological manifold for some neighbourhood $U_1$ of $0\in \mathbb{R}^n$. On the other hand \eqref{eqeq4} gives
$$\det \left[ \frac{\partial^2 f }{\partial x_i \partial x_j}(a_{\nu})\right]\neq 0,\quad \textrm{ dla }\nu\in \mathbb{N},$$
hence and from \eqref{eqeq5} $(f-F)$ has Morse singularities in points $a_{\nu}$, so, $(f-F)^{-1}(0)$ is not $(n-1)$-dimensional topological manifold in any neighbourhood of point $a_{\nu}$. This contradiction completes the proof of Theorem \ref{BLplus1}.

\begin{flushleft}
\scriptsize{
\noindent Piotr Migus\newline 
Faculty of Mathematics and Computer Science, University of \L \'od\'z, 
S. Banacha 22, 90-238 \L \'od\'z, POLAND \newline
 \indent E-mail: migus@math.uni.lodz.pl}
\end{flushleft}

\begin{flushleft}
\scriptsize{
Tomasz Rodak\\
Faculty of Mathematics and Computer Science, University of \L \'od\'z, 
S. Banacha 22, 90-238 \L \'od\'z, POLAND \newline 
\indent E-mail: rodakt@math.uni.lodz.pl}
\end{flushleft}

\begin{flushleft}
\scriptsize{
Stanis\l{}aw Spodzieja\newline 
Faculty of Mathematics and Computer Science, University of \L \'od\'z, 
S. Banacha 22, 90-238 \L \'od\'z, POLAND \newline
 \indent E-mail: spodziej@math.uni.lodz.pl}
\end{flushleft}
\end{document}